\documentclass[11pt,english]{article}
\usepackage[T1]{fontenc}
\usepackage[latin1]{inputenc}

\makeatletter

\newcommand{\binom}[2]{{#1 \choose #2}}

\providecommand{\tabularnewline}{\\}

\title{The Semigroup of a Word}
\author{Peter M. Higgins \& Norman R. Reilly}
\date{}
\def\qed{\quad\vrule height4.17pt width4.17pt depth0pt}

\makeatother

\usepackage{babel}
\begin{document}

\title{Burrows-Wheeler transformations and de Bruijn words}

\author{Peter M. Higgins, University of Essex, U.K. }
\maketitle
\begin{abstract}
We formulate and explain the extended Burrows-Wheeler transform of
Mantaci et al from the viewpoint of permutations on a chain taken
as a union of partial order-preserving mappings. In so doing we establish
a link with syntactic semigroups of languages that are themselves
cyclic semigroups. We apply the extended transform with a view to
generating de Bruijn words through inverting the transform. We also
make use of de Bruijn words to facilitate a proof that the maximum
number of distinct factors of a word of length $n$ has the form $\frac{1}{2}n^{2}-O(n\log n)$.
\end{abstract}

\section{Introduction}

\subsection{Definitions and Example}

The original notion of a Burrows-Wheeler (BW) transform, introduced
in {[}2{]}, has become a major tool in lossless data compression.
It replaces a primitive word $w$ (one that is not a power of some
other word) by another word $BW(w)$ of the same length over the same
alphabet but in a way that is generally rich in letter repetition
and so lends to easy compression. Moreover the transform can be inverted
in linear time; see for example {[}3{]}. Unfortunately, not all words
arise as Burrows-Wheeler transforms of a primitive word so, in the
original format, it was not possible to invert an arbitrary string.
The extended BW transform however does allow the inversion of an arbitrary
word and the result in general is a multiset (a set allowing repeats)
of \emph{necklaces}, which are conjugacy classes of primitive words.
This was first explicitly introduced in {[}8{]} by Mantaci et. al.
based on the bijection between these two collections first enunciated
by Gessel and Reutenauer in {[}5{]}.

In this opening section we will explain and prove the existence of
the extended transform in a fashion that emphasises the approach whereby
a permutation on a finite chain is expressed as a disjoint union of
one-to-one partial order-preserving mappings.

\textbf{Notation and Background }The underlying base set for our mappings
will be the finite chain $[n]=\{0<1<\cdots<n-1\}$. As usual $A^{*}$
will stand for the free monoid over $A=\{a_{0},a_{1},\cdots\}$, which
is simply the set of all words, or strings, over the alphabet $A$
together with the empty word $\varepsilon$, although throughout this
paper we assume a fixed order $a_{0}<a_{1}<\cdots$ for $A$. The
free semigroup is denoted by $A^{+}=A^{*}\setminus\{\varepsilon\}$.
For emphasis, we sometimes denote equality of $u,v\in A^{+}$ by $u\equiv v$.
The set of letters that occur at least once in $w\in A^{*}$ is known
as the \emph{content} of $w$, denoted by $c(w)$. Following {[}8{]}
we shall denote the first and last letters of a word $w\in A^{+}$
respectively by $F(w)$ and $L(w)$. In general, the $i$th letter
of a word $w$ is written as $(w)_{i}$. The number of instances of
the letter $a_{i}$ in a word $w$ will be denoted by $|w|_{a_{i}}$,
while the length of $w$ is written $|w|$. We say that $w$ is \emph{primitive}
if $w$ is not a power of some other word. A word $u\in A^{+}$ is
a \emph{factor }of $w\in A^{+}$ if $w\in A^{*}uA^{*}$; $u$ is an
\emph{$m$-factor} of $w$ if additionally $u\in A^{m}$. We call
$u$ a \emph{prefix }(respectively \emph{suffix})\emph{ }of $w$ if
$w\in uA^{*}$ (respectively $w\in A^{*}u$). A \emph{subword} of
$w$ is any word that may be formed by deletion of some of the letters
of $w$; it follows that the factors of $w$ represent a special class
of subwords of $w$.

A standard text for results concerning combinatorics on words is {[}7{]}
in which may be found proofs for simple unproved assertions concerning
roots and conjugates that follow. If $w=uv,(u,v\in A^{*})$ we say
that $w'=vu$ is a \emph{conjugate} of $w$. The relation $\sim$
on $A^{*}$ whereby $w\sim w'$ if $w'$ is a conjugate of $w$ is
an equivalence relation on $A^{*}$. In the case of a primitive word
$w$, the equivalence classes of $\sim$ are known as \emph{necklaces,
}and we denote the necklace of a word $w$ by $n(w)$; the \emph{length
}of $n(w)$ is $|w|$, which is also the cardinal of the necklace
as $w$ is primitive. The first word of $n(w)$ in the lexicographic
order is known as its \emph{Lyndon word.} A \emph{border} of a word
$w$ is word $u\in A^{+}$ such that $w\in uA^{+}\cap A^{+}u$. No
Lyndon word has a border (see Proposition 2.2(iii)). 

The \emph{root} of a word $w$ is the shortest factor $r=$ root$(w)$
of $w$ such that $w=r^{t}$ for some $t\geq1$. Two words $w$ and
$u$ commute in $A^{+}$ if and only if they share a common root,
which is in turn equivalent to the condition that $w$ and $u$ have
a common power. The number of distinct conjugates of a word $w$ equals
the length of root$(w)$ and root$(w')$, the root of a conjugate
$w'$ of $w$, is a conjugate of root$(w)$.

For a word $w$ we denote the infinite one-sided word $www\cdots$
by $w^{\omega}$ with the notion of factor extending in the obvious
way. Note that $u^{\omega}=v^{\omega}$ if and only if root$(u)=$
root$(v)$. The factors $u$ of $w^{\omega}$ of finite length are
the \emph{power factors }of $w$; a power factor for which $|u|\leq|w|$
is a \emph{cyclic factor} of $w$: equivalently $u$ is a factor of
some conjugate of $w$. 

The interval $I=[i,j]$ of a chain $X$ is the subset $I=\{k:i\leq k\leq j\}$.
A mapping $\alpha$, the domain and range of which are both subsets
of $A$, is \emph{order-preserving} if when $a\cdot\alpha$ and $b\cdot\alpha$
are both defined, $\alpha$ satisfies the condition: 
\[
a\leq b\rightarrow a\cdot\alpha\leq b\cdot\alpha\:(a,b\in A).
\]
 We shall frequently use the action notation, $a\cdot\alpha$ as opposed
to juxtaposition $a\alpha$ when the symbol on the right is a function
and not a product in $A^{*}$ (although a central dot is also used
at times simply as a visual separator within a word). Mapping composition
is written from left to right. Here we write $PI_{n}$ to denote the
(inverse) semigroup of all partial one-to-one mappings on $[n]$,
and we denote the (inverse) subsemigroup of all order-preserving members
of $PI_{n}$ by $POI_{n}$. 

\textbf{Example 1.1 }We give a example, following {[}8{]}, that illustrates
how to effect the bijection from multisets of necklaces to words and
how to reverse this process. Let our alphabet be $A=\{a<b\}$ and
let the set of Lyndon words of our necklaces be $M=\{aab,\,ab,\,abb\}$.
Consider the collection of all words of the form $u^{\frac{l}{|u|}}$,
where $u\in n(v)$ $(v\in M)$ and $l$ is the least common multiple
of the lengths of the words of $M$: in this instance $l=3\times2=6$.
All these words then have common length $l$. We order this set of
words lexicographically to yield, in our example, the following array.

{\center%
\begin{tabular}{|c|c|c|c|c|c|}
\hline 
a & a & \textbf{b} & a & a & \textbf{b}\tabularnewline
\hline 
\hline 
a & b & \textbf{a} & a & b & \textbf{a}\tabularnewline
\hline 
a & \textbf{b} & a & b & a & \textbf{b}\tabularnewline
\hline 
a & b & \textbf{b} & a & b & \textbf{b}\tabularnewline
\hline 
b & a & \textbf{a} & b & a & \textbf{a}\tabularnewline
\hline 
b & \textbf{a} & b & a & b & \textbf{a}\tabularnewline
\hline 
b & a & \textbf{b} & b & a & \textbf{b}\tabularnewline
\hline 
b & b & \textbf{a} & b & b & \textbf{a}\tabularnewline
\hline 
\end{tabular}\endcenter}

The Burrows-Wheeler transform of $M$ is then the word formed by the
$l$th column of the table, read from the top, which in this case
gives $BW(M)=babbaaba$. The word $BW(M)$ is also formed by the list
of last letters $L(u)$: both renditions of $BW(M)$ are highlighted
in bold in the table. In {[}8{]} $BW(M)$ was defined by the letters
$L(u)$. Their definition was also framed in context of the infinite
table $T$ of rows $u^{\omega}$, which simply consists of the table
of the first $l$ columns of $T$, as defined above, repeated infinitely
often. However, as explained in {[}8{]}, the table does not need to
be extended to $l$ columns in order to determine the order of the
rows: by a theorem of Wilf and Fine on word periodicity, the order
of two rows that are respective powers of the root words $u$ and
$v$ matches the lexicographic order of their prefixes of length $k=|u|+|v|-\mbox{gcd\ensuremath{(|u|,|v|)}}$
(and this bound is tight). Hence the number of columns required in
order to determine the row order of the table is always less than
the sum of the lengths of the longest two necklaces of the multiset.
The formal use of the lcm $l$ here allows us to define $BW(M)$ as
a specified column of the table, which is a conceptual convenience
used in our proofs. The stipulation that the words of $M$ be primitive
is necessary in order that the $BW$ transform be one-to-one. Note
that the roots of the words are not in lexicographic order: the root
$baa$ precedes the root $ba$ in the table. However the Lyndon roots
do appear in lexicographic order: $aab<ab<abb$ both lexicographically
and in the rows of the table (see Theorem 1.2.13).

We recover the set $M$ from $w=BW(M)$ by way of the so-called \emph{standard
permutation }$\pi=\pi(w)$. To construct $\pi$, take the first column
of the table, which consists of the content of the words of $M$ arranged
in alphabetical order with the number of occurrences of a letter equal
to the number of instances of that letter among the Lyndon words of
$M$. In our example the column of first letters forms the word $F(M)=aaaabbbb$.
The permutation, $\pi(w)$ is then the union of a collection of partial
one-to-one and order-preserving mappings, one for each member of $c(w)$.
In this case $\pi=\pi_{a}\cup\pi_{b}$; the domain and range of $\pi_{a}$
is defined respectively by the positions of the instances of the letter
$a$ in $F(M)$ and $BW(M)$ respectively. Since $\pi_{a}$ is one-to-one
and order-preserving, $\pi_{a}$ is defined uniquely by its domain
and range, and of course $\pi_{b}$ is defined in the same fashion,
and so on for any remaining letters in $c(w).$ In our example we
obtain: 
\[
\pi(w)=\binom{0\,1\,2\,3\,4\,5\,6\,7}{1\,4\,5\,7\,0\,2\,3\,6}=(0\,1\,4)(2\,5)(3\,7\,6),
\]
with dom $\pi_{a}=\{0,1,2,3\}$ and dom $\pi_{b}=\{4,5,6,7\}$. The
cardinality of $M$ is equal to the number of cycles in the disjoint
cycle representation of $\pi$, which here is $3$. We may retrieve
the Lyndon word of the multiset $M$ corresponding to each cycle of
$\pi(w)$ by simply replacing each integer $m$ in the cycle by the
letter $c\in A$ such that $m\in$ dom $\pi_{c}$. In our case this
means that we write $a$ whenever we see a number from $0$ to $3$
and we write $b$ otherwise. In this way we recover $M=\{aab,\,ab,abb\}$.

\subsection{Establishing the transform through partial order-preserving mappings}

Using Example 1.1 as a guide, we formally define the Burrows-Wheeler
transform and explain its inversion. 

\textbf{Definition 1.2.1} \emph{(Conjugation Map)} Let $\Pi:A^{+}\rightarrow A^{+}$
be the mapping whereby $au\mapsto ua$ $(a\in A,u\in A^{*})$.

\textbf{Proposition 1.2.2} The Conjugation Map has the following properties:

(i) $\Pi$ is a permutation on $A^{+}$; 

(ii) if $S\subseteq A^{*}$ is closed under conjugation then $\Pi|_{S}$
permutes $S$.

(iii) Suppose that $S\subseteq aA^{n}$$(a\in A,n\geq0)$. Then $\Pi$
acts in an order-preserving manner on $S$.

(iv) For any word $w$ with root$(w)=r$, $|r|$ is the least positive
integer $t$ such that $w\cdot\Pi^{t}=w$.

\emph{Proof} (i) is clear from the definition and (ii) follows from
(i) as the given condition ensures that $S$ is closed under both
$\Pi$ and $\Pi^{-1}$. To see (iii) suppose that $au\leq av$ with
$au,av\in aA^{n}$. Since $|u|=|v|$, it follows that $u\leq v$ whence
$ua\leq va$ and so $\Pi$ is order-preserving on the set $aA^{n}$.
As for (iv), if $w\equiv xy$ then $w\cdot\Pi^{|x|}=yx$ so in particular
$w\cdot\Pi^{|r|}=w$. Suppose that $1\leq|x|<|r|$ so that $r\equiv xx'$
say. Then $w'=w\cdot\Pi^{|x|}\in x'xA^{*}$ and since $|x'x|=|r|$
but $x'x\not\equiv r$ as $r$ is primitive, it follows that $w'\neq w$.
$\qed$

\textbf{Definitions 1.2.3 }\emph{(Burrows-Wheeler map)} Let ${\cal M}$
denote the set of all finite multisets of necklaces over $A$. Let
$BW:{\cal M}$$\rightarrow A^{*}$ denote the Burrows-Wheeler map,
the action of which is defined as follows. Take any $M\in{\cal M}$
so that $M=\{n_{1},n_{2},\cdots,n_{t}\}$ $(t\geq0)$ and let $l$
be the least common multiple of the lengths of the $n_{i}$. Sort
by lexicographic order the collection $T=T(M)$ of powers $u^{\frac{l}{|u|}}$,
where $u$ is a word of the necklace $n_{i}.$ The table $T$ is then
a dictionary of $n=|n_{1}|+|n_{2}|+\cdots+|n_{t}|$ words of common
length $l$. The word $BW(M)$ is then the final column, read from
top to bottom, of $T$. (Conventionally, $BW$ maps the empty set
to the empty word.)

\textbf{Definition 1.2.4 }\emph{(Standard permutation of a word) }Let
$w\in A^{n}$ and let $f(w)$ be the rearrangement of the letters
of $w$ in lexicographic order. For each letter $a\in c(w)$ we define
a partial one-to-one order-preserving mapping $\pi_{a}\in PIO_{n}$
through specifying dom $\pi_{a}$ and ran $\pi_{a}$ as follows: dom
$\pi_{a}$ is the interval of length $|w|_{a}$ corresponding to the
positions occupied by $a$ in $f(w)$ while ran $\pi_{a}$ is the
set of positions occupied by $a$ in $w$. The \emph{standard permutation}
of $w$ is then $\pi=\cup_{a\in c(w)}\pi_{a}$. 

\textbf{Remark 1.2.5 }For any $i\in[n]$ there is a unique $a\in A$
such that $i\cdot\pi=i\cdot\pi_{a}$. For any $u\in A^{*}$, $u\equiv b_{1}b_{2}\cdots b_{m}$
we may define $\pi_{u}=\pi_{b_{1}}\pi_{b_{2}}\cdots\pi_{b_{m}}$.
We note that $\pi_{u}\in PIO_{n}$ and for any $m\geq1$ and $i\in[n]$
there is a unique word $u=u_{i,m}$ of length $m$ such that $i\cdot\pi_{u}$
is defined. 

\textbf{Proposition 1.2.6 }{[}9, Proposition 10{]} Let $M\in{\cal M}$
as in Definition 1.2.3, let the set of words that form the rows of
$T(M)$ be denoted by $R(M)$ and let $u_{i}\in R(M)$ $(i\in[n]$).
Let $\pi=\pi(w)$ be the standard permutation of $w=BW(M)$. Then
the mapping $u_{i}\mapsto u_{i\cdot\pi}$ is the restriction of the
conjugation map $\Pi$ to $R(M)$. 

\emph{Proof} Suppose that $F(u_{i})=a$ and that $u_{i}$ is the $j$th
word of $R(M)\cap aA^{*}$. Then the $j$th instance of $a$ in the
first column of $T(M)$ occurs in row $i$. Hence, regarded as intervals
of $[n]$, dom $\Pi|_{R(M)\cap aA^{*}}=$ dom $\pi_{a}$. Similarly,
since $w$ is the final column of $T(M)$, ran $\pi_{a}=R(M)\cap A^{*}a=(R(M)\cap aA^{*})\Pi$.
Therefore since $\pi_{a}$ and $\Pi|_{R(M)\cap aA^{*}}$ are order-preserving
mappings (the latter by Proposition 1.2.2(iii)) with common domain
and range, they are equal. Since this is true for all letters $a\in A$,
we infer that $\pi=\Pi|_{R(M)}$ in that $i\mapsto i\cdot\pi$ if
and only if $u_{i}\mapsto u_{i\cdot\pi}$ under $\Pi$. $\qed$

The following was observed in {[}3{]}, at least for the case of the
BW transform of a single necklace.

\textbf{Proposition 1.2.7 }Let $w=BW(M)\in A^{n}$, let $\pi=\pi(w)$
be the standard permutation and let $T(M)=(a_{ij})$. Then $a_{ij}=a_{i\cdot\pi,j-1}$,
which is to say that $\pi$ maps each column of $T(M)$ to its predecessor
column modulo $l$, the number of columns of $T(M)$. In particular
$\pi$ maps the first column of $T(M)$ to the last.

\emph{Proof} Let $u_{i}$ be a row of $T(M)$ with $F(u_{i})=b$ so
that $u_{i}=bu$ say. Then by Proposition 1.2.6, $u_{i\cdot\pi}=ub$.
The letter $a=a_{ij}$ will therefore be shifted one place back to
appear in column $j-1$ and in row $i\cdot\pi$ so that $a=a_{ij}=a_{i\cdot\pi,j-1}$.
$\qed$

\textbf{Definition 1.2.8 }\emph{(Table of a word) }Let $w=b_{0}b_{1}\cdots b_{n-1}\in A^{+}$
and let $\pi=\pi(w)$ be the standard permutation of $w$. Let us
write the cycle $C_{i}=(i\,\,i\cdot\pi\,\,i\cdot\pi^{2}\,\,\cdots\,\,i\cdot\pi^{r-1})$
so $r$ is least such that $i\cdot\pi^{r}=i$ and let $l$ denote
the lcm of the cycle lengths. Define the table $T(w)$ to be the $n\times l$
table, the $i$th row of which is the unique word $u=u_{i}\in A^{l}$
such that $i\cdot\pi_{u}$ is defined.

\textbf{Proposition 1.2.9 }Let $w,\pi$ and $T(w)$ be as in Definition
1.2.8. Let $r=r(i)$ be the length of $C_{i}$ and let $x\in A^{r}$
be the corresponding prefix of $u=u(i)$, the $i$th row of $T(w)$.
Then

(i) $x$ is the root of $u$;

(ii) all conjugates of $x$ arise as roots of the rows of $T(w)$
with multiplicity equal to that of $x$.

(iii) The rows of $T(w)$ are ranked lexicographically.

(iv) The final column of $T(w)$ is $w$.

\emph{Proof} (i) By construction, $i\cdot\pi^{r}=i\cdot\pi_{x}=i$
and $x$ is the shortest prefix of $u$ with this property. In particular,
it follows from this that $u=x^{\frac{l}{r}}$. To show that $x$
is itself primitive, and so the root of $u$, suppose to the contrary
that $x=y^{t}$ for some $t\geq2$. Then $i\cdot\pi_{y}\neq i$; without
loss suppose that $i<i\cdot\pi_{y}$. By applying $\pi_{y}$ to both
sides of this inequality (remembering that $i\cdot\pi_{y^{s}}$ is
defined for all $s\leq t)$ we infer that 
\[
i<i\cdot\pi_{y}<i\cdot\pi_{y^{2}}<\cdots<i\cdot\pi_{y^{t}}=i\cdot\pi_{x}=i,
\]
 a contradiction. Hence $t=1$ and $x$ is the root of $u$, as claimed. 

(ii) Let $y=qp$ be a conjugate of $x=pq$, the root of $u(i)$. Then
\[
(i\cdot\pi_{p})\cdot\pi_{y}=i\cdot\pi_{py}=i\cdot\pi_{xp}=i\cdot\pi_{x}\pi_{p}=i\cdot\pi_{p}
\]
and since $y$ is primitive, it follows that $u(i\cdot\pi_{p})=y^{\frac{l}{|y|}}$
and $y$ is indeed the root of $u(i\cdot\pi_{p})$. This process associates
each instance of the root $x$ with an instance of the conjugate $y$
in a one-to-one fashion, thereby matching the multiplicity of $x$
to that of each of its conjugates $y$ in the table $T(w).$

(iii) Let $i<j$, let $u=u(i)$ and $v=v(j)$ be distinct words that
occupy the respective rows $i$ and $j$ of $T(w)$ and let $p\in A^{*}$
be the longest common prefix of $u$ and $v$ so that $u=pu_{1}$
and $v=pv_{1}$ say. Then since $\pi_{p}$ is order-preserving we
have $i_{1}=i\cdot\pi_{p}<j\cdot\pi_{p}=j_{1}$. Since $u$ and $v$
have common length $l$, it follows that $F(u_{1})=a,\,F(v_{1})=b$
say with $a\neq b$. Moreover, since $i_{1}\in$ dom $\pi_{a},\,j_{1}\in$
dom $\pi_{b}$ and $i_{1}<j_{1}$, it follows that $a<b$ and so $u<v$,
as required. 

(iv) Let $(a_{ij})$ denote the table $T(w)$. Then $a_{ij}=a$ if
and only if $i\cdot\pi^{j-1}\in$ dom $\pi_{a}$. In particular, taking
$j=l$ gives that $i\cdot\pi^{l-1}\in$ dom $\pi_{a}$, whence $i\cdot\pi^{-1}\in$
dom $\pi_{a}$. At the same time we observe that $(w)_{i}=a$ exactly
when $i\pi^{-1}\in$ dom $\pi_{a}$ and therefore $a_{il}=(w)_{i}$
for all $i\in[n]$, whence $w$ is indeed the final column of $T(w)$.
$\qed$

\textbf{Definition 1.2.10 }\emph{(Inverse Burrows-Wheeler map) }Define
$I:A^{*}\rightarrow{\cal M}$ as follows. Given $w\in A^{n}$, form
$T(w)$ as in Definition 1.2.8. Let $M=I(w)$ be the set of necklaces
defined by the roots of the rows of $T(w)$. (With $\varepsilon\mapsto\emptyset$
under $I$.)

\textbf{Theorem 1.2.11} {[}5, 8{]} The mapping $I$ of Definition
1.2.10 is the inverse Burrows-Wheeler transform $BW^{-1}:A^{*}\rightarrow{\cal M}$. 

\emph{Proof }We first prove that for any $M\in{\cal M}$, $I(BW(M))=M$.
Let $T=T(M)$ be the table of $M$ and let $w=BW(M)\in A^{n}$ as
in Definition 1.2.3. We show that the $i$th row $u=u_{i}$ of $T(M)$
is the $i$th row of $T(w)$. By Proposition 1.2.6, identifying the
rows of $T(M)$ with the chain $[n]$ allows us to say that $\pi(w)=\Pi|_{R(M)}$.
In particular the lcm of the cycle lengths of both permutations is
a common value $l$, and by Proposition 1.2.2(iv) $l$ is the lcm
of the lengths of the roots of the words of $R(M),$ so that $T$
is an $n\times l$ array.

Now suppose that $u=av$ $(a\in A)$. By Proposition 1.2.6 it follows
that $va=u_{i\cdot\pi}=u_{i\cdot\pi_{a}}$, so that $i\cdot\pi=i\cdot\pi_{a}$.
Repeated application of this observation gives that $i\cdot\pi^{l}=i\cdot\pi_{u}$
so that $u$ is the unique word of length $l$ such that $i\cdot\pi_{u}$
is defined. Hence $T(M)=T(w)=T$ say. By Definition 1.2.10, $I(w)$
is the set of necklaces formed by the roots of $T$, which is the
set $M$ itself, and so $I(BW(M))=M$. 

Conversely, take any $w\in A^{n}$ say and let $M=I(w)$. By Definition
1.2.10, $M$ is the collection of necklaces of the roots of the rows
of $T(w)$. By Proposition 1.2.9(i), if $x$ is the root of row $i$
in $T(w)$, then $r=|x|$ is the length of the cycle $C_{i}$ of $\pi(w)$.
It follows that there is a common value $l$ for the lcm of the lengths
of the roots of the rows of $T(w)$ (which is the row length of $T(M)$)
and the lcm of the cycle lengths of $\pi(w)$ (which is the row length
of $T(w)$). By Proposition 1.2.9(ii), all members of $n(x)$ appear
as roots of rows of $T(w)$ with equal multiplicity while by (iii)
the rows of $T(w)$ are ranked lexicographically. It follows from
all this that $T(w)=T(M)=T$ is an $n\times l$ array. Now $BW(M)$
is the final column of $T$, which by Proposition 1.2.9(iv) is the
word $w$. We conclude that $BW(I(w))=w$. $\qed$ 

\textbf{Remark 1.2.12 T}he first part of the previous proof establishes
that $T(w(M))=T(M)$ while the third paragraph shows that $T(M(w))=T(w)$
so that the bijection between words and necklaces is through equality
of the corresponding table $T$. Moreover Proposition 1.2.6 shows
that the action of $\Pi$ on $R(T)$ corresponds to that of $\pi(w)$
on $[n]$ and Proposition 1.2.7 shows that $\pi$ acts to map each
column of $T$ onto its predecessor modulo $l$. 

\textbf{Theorem 1.2.13 }Let $M\in{\cal M}$, $T=T(M)$ and let $i<j$
with $u=u(i),v=u(j)$ two words in the set of rows $R(M)$ of $T$.
Then $u<r=$ root$(v)$ if root$(u)$ is Lyndon. In particular the
Lyndon words appear in $R(T)$ in lexicographic order.

\emph{Proof} We prove the first statement by showing that if $r\leq u(i)$
then root$(u)$ is not Lyndon. Given this claim, suppose that root$(u)$
and root$(v)$ are both Lyndon words such that $u<v$. Then root$(u)\leq u<$
root$(v)$ so that the Lyndon roots do indeed appear in lexicographic
order in $T$.

Since $u<v$ with $r=$ root$(v)\leq u$ it follows that $v$ is not
primitive and so $v=r^{t}$ for some $t\geq2$. Since $|u|=|v|$ and
$u<v$ we may write $u=pax,v=pby$ with $a,b\in A,\,p,x,y\in A^{*}$
and $a<b$. If $|p|<|r|$, then $r=pbq$ say whence $u<r$, contrary
to hypothesis and so $|r|\leq|p|$ whence, since $v$ is a power of
$r$, $p=r^{m}s$ for some maximal $m\geq1$, and where $s\in A^{*}$
is a prefix of $r$. It follows that $r=st$ where $F(t)=b$ so that
$t=bw$ say $(w\in A^{*})$ whence $r=sbw$. Taking the factorization
$u=r^{m}sax$, we see that $u'=u\cdot r^{m}=sa(xr^{m})$ is a conjugate
of $u$. We also have the factorization $u=r^{m}sax=sb(wr^{m-1}sax$),
whence $u'<u$ as $sa<sb$, which implies that root$(u')<$ root$(u)$
and so root$(u)$ is not Lyndon, as required. $\qed$

\section{Semigroup of the Burrows-Wheeler transform}

\subsection*{Semigroup of a necklace}

In {[}6{]} the author wrote about the semigroup $S(u)$ generated
by the letters acting by conjugation on the necklace of a primitive
word $u$. In particular the question of when two words $u$ and $v$
have isomorphic semigroups $S(u)$ and $S(v)$ was settled by Theorem
2.4 of {[}6{]}. The semigroup $S(u)$ is exactly the semigroup generated
by the partial mappings $\pi_{a}$ $(a\in c(u))$ encountered above.
We show here that $S(u)$ is isomorphic to the syntactic semigroup
of the cyclic semigroup generated by the word $u$. 

We begin with a fixed primitive word $u\in A^{n}$ over the finite
ordered alphabet $A=\{a_{0}<a_{1}<\cdots<a_{k-1}\}$. Consider the
necklace $n(u)=\{u_{0}<u_{1}<\cdots<u_{n-1}\}$, ordered lexicographically. 

\textbf{Definition 2.1 }Identify the chain $n(u)$ with the chain
$[n]$. The semigroup $S(u)$ is the subsemigroup of $POI_{n}$ generated
by the set of $k$ partial mappings $\{\pi_{a_{i}}\}$ where $\pi_{BW(n(u))}=\cup_{i=0}^{k-1}\pi_{a_{i}}$
$(a_{i}\in A)$. 

In this section it is convenient to denote the mapping $\pi_{a}$
by $a'$ so that the semigroup $S(u)$ is generated by the set of
partial mappings $a'$ ($a\in A)$ where $u_{j}\in$ dom $a'$ if
and only if $F(u_{j})=a$ so that $u_{j}=ax$ say in which case $(ax)a'=xa\in n(u)$.
We write this using action notation as $ax\cdot a=xa$, allowing us
to suppress the dash to the right of the central dot without introducing
ambiguity. The free monoid $A^{*}$ acts on the right of $n(u)$ in
that $u_{j}\cdot(xy)=(u_{j}\cdot x)\cdot y$ for all $u_{j}\in n(u)$
and $x,y\in A^{*}$ (taking $\varepsilon'$ to be the identity mapping).
Note that $S(u)$ depends only on the necklace $n(u)$ and not its
representative (and so we may assume that $u=u_{0},$ the Lyndon word
of $n(u),$ although this is not necessary). We make use of the following
facts from Proposition 1.3 in {[}6{]}; part (iii) is well-known -
see for example the text {[}7{]}.

\textbf{Proposition 2.2 }Let $u=b_{1}b_{2}\cdots b_{n}\in A^{+}$
and $t\geq0$ be an integer. Let $z=u^{m}b_{1}b_{2}\cdots b_{s}$
be the prefix of $u^{\omega}$ of length $t$ so that $t=mn+s$$(0\leq m,\,0\leq s\leq n-1).$
Write $v=b_{1}b_{2}\cdots b_{s}$ and define $w\in A^{+}$ by $u=vw$.
Then

(i) $u\cdot v=wv$ and $u\cdot u=u$;

(ii) $z\equiv u^{m}v$ is the unique word $y$ of length $t$ such
that $u\cdot y$ is defined.

(iii) A Lyndon word $u$ has no border.

\emph{Proof} (i) and (ii) are immediate consequences of the definition
of the action of each letter on a given word. As for (iii), suppose
to the contrary that $u=xv\equiv vw$ for some $x,v,w\in A^{+}$.
Then since $u$ is Lyndon (and primitive) we may apply (i) to infer
that $u<u\cdot x$ and $u<u\cdot v$. From the first of these inequalities
we get $u\cdot v<u\cdot xv$ as the latter is defined because $u\cdot xv=u\cdot u=u$.
However we then obtain $u<u\cdot v<u\cdot xv=u\cdot u=u$, which is
a contradiction. Therefore $u$ has no border. $\qed$

We now introduce a second realisation of $S(u)$ via a certain syntactic
congruence, thus producing $S(u)$ without reference to mappings.
(For background on syntactic semigroups and congruences see {[}11{]}.)
Let $\langle u\rangle$ be the subsemigroup of $A^{+}$ of all positive
powers of $u$. Let $\rho=\rho_{u}$ be the \emph{syntactic congruence}
on $A^{+}$ generated by $\langle u\rangle$ so that for $x,y\in A^{+}$:
\begin{equation}
x\rho y\leftrightarrow(pxq\in\langle u\rangle\leftrightarrow pyq\in\langle u\rangle\,\forall p,q\in A^{*})
\end{equation}

\textbf{Definition 2.3 }The semigroup $S_{u}=A^{+}/\rho_{u}$.

\textbf{Lemma 2.4 }Let $u\equiv wv$ and $u'\equiv vw$ be conjugate
words. Then $S_{u}=S_{u'}$. 

\emph{Proof}\textbf{\emph{ }}Suppose that $(x,y)\in\rho_{u}$. Then
for any $p,q\in A^{*}$ we have that if $pxq\equiv u'^{m}\equiv(vw)^{m}$
for some $m\geq1$ then$(wp)x(qv)\equiv w(vw)^{m}v\equiv(wv)^{m+1}\equiv u^{m+1}$.
Since $(x,y)\in\rho_{u}$ this in turn implies that $(wp)y(qv)\equiv u^{r+1}\equiv(wv)^{r+1}$
for some $r\geq1$ $(r\neq0$ as $y\neq\varepsilon)$, whence $pyq\equiv(vw)^{r}\equiv u'^{r}$.
Hence it follows that $pxq\in\langle u'\rangle$ implies that $pyq\in\langle u'\rangle$.
Interchanging the roles of $x$ and $y$ in this argument yields the
conclusion that $\rho_{u}\subseteq\rho_{u'}$ and by symmetry of the
conjugation relation we see that the reverse inclusion also holds.
Therefore $\rho_{u}=\rho_{u'}$ and $S_{u}=S_{u'}$. $\qed$

\textbf{Theorem 2.5 }For any primitive word $u$, $S_{u}\cong S(u).$ 

\emph{Proof }For each $x\in A^{+}$, let $[x]=x\rho$ be the corresponding
member of $S_{u}$ and $x'$ be that of $S(u)$. We show that a required
isomorphism is given by the mapping $\theta:[x]\mapsto x'$. We first
verify that $[x]=[y]$ if and only if $x'=y'$, thereby showing that
$\theta$ is an injective function. It is then clear from the definition
that $\theta$ is also surjective and $\theta$ is a homomorphism
as for any $x,y\in A^{+}$ we then have 
\[
([x][y])\theta=[xy]\theta=(xy)'=x'y'=[x]\theta[y]\theta.
\]
To this end suppose that $x\rho y$ and suppose further that $vw\in n(u)$,
where $u=wv$ and that $vw\cdot x$ is defined. By Proposition 2.2(ii),
$x\equiv(vw)^{m}c$ for some $m\geq0$, where $vw\equiv cd$ say $(c,d\in A^{*})$.
We shall show that $(vw)\cdot x=(vw)\cdot y$: 
\begin{equation}
vw\cdot x=(vw)\cdot(vw)^{m}c=(vw)\cdot c=(cd)\cdot c=dc\in n(u)
\end{equation}
where the second and fourth equalities are by Proposition 2.2(i).
Then since $vw\equiv cd$ we have: 
\[
wxdv\equiv w(vw)^{m}cdv\equiv(wv)^{m}w(vw)v\equiv(wv)^{m+2}
\]
and so $wxdv\in\langle u\rangle$. Therefore since $x\rho y$ we infer
that $wydv\equiv(wv)^{r}$ for some $r\geq2$ ($r>1$ as $y\neq\varepsilon$.)
Hence, by cancelling $w$ on the left and $v$ on the right of this
equation we obtain: 
\begin{equation}
yd\equiv(vw)^{r-1}
\end{equation}
Invoking (2) and then (3) we infer that 
\begin{equation}
vw\cdot xd=(vw\cdot x)\cdot d=dc\cdot d=cd;\,\,\,\,vw\cdot yd=vw\cdot(vw)^{r-1}=vw\equiv cd
\end{equation}
Since the mapping $d'$ is injective, (4) allows us to deduce that
$vw\cdot x=vw\cdot y$. Since $x,y\in A^{+}$ were arbitrary, it follows
that $x\rho y$ implies that $x'=y'$ as the argument shows that for
any $u_{i}=vw\in n(u)$, if one of $u_{i}\cdot x,\,u_{i}\cdot y$
is defined, then both are defined and are equal.

To prove the converse we next suppose that for some $x,y\in A^{+},$
$x'=y'$ and suppose further that $pxq\equiv u^{m}$ for some $p,q\in A^{*}$
and $m\geq1$. We verify that $pyq\in\langle u\rangle$. The following
argument will hold with the roles of $x$ and $y$ reversed and so
this claim yields that if $x'=y'$ then $x\rho y$, thus establishing
that $\theta$ is a one-to-one mapping from $S_{u}$ into $S(u)$.
Since $x'=y'$ we obtain$(u\cdot p)\cdot x=(u\cdot p)\cdot y\Rightarrow((u\cdot p)\cdot x)\cdot q=((u\cdot p)\cdot y)\cdot q$$\Rightarrow u\cdot(pxq)=u\cdot(pyq)\Rightarrow u\cdot u^{m}=u\cdot(pyq);$
by Proposition 2.2(i) we infer that $u=u\cdot pyq$. By Proposition
2.2(ii), $pyq\equiv u{}^{s}v$ $(s\geq0)$ for some non-empty prefix
$v$ of $u\equiv vw$ say. However then we obtain 
\[
u\cdot pyq=u\cdot u^{s}v=u\cdot v=vw\cdot v=wv;\,\,u\cdot pyq=u\equiv vw.
\]
Hence $u=wv\equiv vw$ and since $u$ is primitive it follows that
$v\equiv u$, $w\equiv\varepsilon$ and so $pyq\equiv u^{s+1}$ for
some $s\geq0$. In particular, $pyq\in\langle u\rangle$, as required
to complete the proof of the claim. Therefore $\theta$ is an isomorphism
from $S_{u}$ to $S(u)$. $\qed$

For a multiset of necklaces $M$, we may define the semigroup $S(M)$
in terms of the partial mappings of the standard permutation of $BW(M)$. 

\textbf{Theorem 2.6 }Let $M=\{n_{i}=n(u_{i})$\} be a multiset of
necklaces and let $n=|n_{1}|+|n_{2}|+\cdots+|n_{t}|$. Let $S(M)$
be the subsemigroup of $POI_{n}$ generated by the set of mappings
$\{\pi_{a}\}$ of $\pi=\pi(BW(M))$. Then $S(M)$ is a subsemigroup
of $POI_{n}$ isomorphic to a subdirect product of the syntactic semigroups
$S_{u_{i}}$.

\emph{Proof}\textbf{ }Let $C$ denote any member of the set of domains
$\{C_{1},C_{2},\cdots,C_{t}\}$ of disjoint cycles of $\pi$. Since
$C\pi=C$ and each $\pi_{a}$ is a restriction of $\pi$, it follows
that $\pi_{a}|_{C}$ is a (possibly empty) one-to-one and order-preserving
mapping in $POI_{C}$, where $C$ inherits a linear order as a subchain
of $[n]$. The mapping whereby $\pi_{a}\mapsto(\pi_{a}|_{C_{1}},\pi_{a}|_{C_{2}},\cdots,\pi_{a}|_{C_{t}})$
induces an injective homomorphism $\phi:S(M)\rightarrow\Pi=POI_{C_{1}}\times POI_{C_{2}}\times\cdots\times POI_{C_{t}}$.
Let $p_{j}$ denote the $j$th projection mapping on $\Pi$ so that
$\phi p_{j}:S(M)\rightarrow POI_{C_{j}}$. We see that $S(u_{i})$
is the image of $\phi p_{j(i)}:S(M)\rightarrow POI_{C_{j(i)}}$ (where
$i\in$ dom $C_{j}$) with generators $\pi_{a}|_{C_{j(i)}}$ $(a\in c(u_{i}))$.
It follows that $\phi$ may be regarded as an injective homomorphism
of $S(M)$ into $S(u_{1})\times S(u_{2})\times\cdots\times S(u_{t})$.
Finally, by Theorem 2.5, $S(u_{i})\cong S_{u_{i}}$, the syntactic
semigroup of $\langle u_{i}\rangle$ and so we conclude that $S(M)$
is isomorphic to a subdirect product of the syntactic semigroups of
each of the languages $\langle u_{i}\rangle$, as required. $\qed$

\section{de Bruijn Words}

In this section we take our alphabet to be $A=\{0<1<\cdots<k-1\}$
$(k\geq2)$, although we continue to refer to its members $a\in A$
as \emph{letters. }An interesting special case is where we take the
BW transform of (the necklace of) a \emph{de Bruijn word of span}
\emph{$n$} over a finite $k$-ary alphabet, which can be defined
as a word $w$ of length $k^{n}$ for which every word of length $n$
appears exactly once as a cyclic factor of $w$. For every $n$ and
for every $k$-ary alphabet $A$, de Bruijn words $d_{n}$ exist and
their number is $\frac{(k!)^{k^{n-1}}}{k^{n}}$ {[}1{]}. 

\textbf{Definition 3.1 }A multiset $M$ of necklaces $\{n_{i}\}$
is a \emph{de Bruijn} \emph{set of span }\textbf{\emph{$n$ }}over
$A$ if $|n_{1}|+|n_{2}|+\cdots+|n_{t}|=k^{n}$ and every $w\in A^{n}$
is a prefix of some power of some word of the necklaces $n_{i}$.

\textbf{Remarks 3.2} The number of distinct prefixes of length $n$
of powers of the words of the necklaces $n_{i}$ is at most $k^{n}$
so, given that $M$ is a de Bruijn set of span $n$, every word in
$A^{n}$ can be read exactly once within the necklaces of $M$. It
also follows in particular that no two necklaces in $M$ are equal
so that $M$ is indeed a set, as opposed to a multiset, of necklaces. 

\textbf{Lemma 3.3 }Let $M$ be a de Bruijn set of span $n$. Then
$M$ contains a necklace of length at least $n$.

\emph{Proof} There exist Lyndon words $u$ of length $n$ (eg. take
$u=ab^{n-1}$, where $a<b$). Let $n_{i}\in M$ be a necklace of cardinal
$m<n$ so that $n=tm+l$ say with $0\leq l\leq m-1$. Any prefix of
length $n$ of a power factor of a word $v\in n_{i}$ has a border
of length $l$ if $l\neq0$ and has a border of length $m$ otherwise.
Since $u$ is a Lyndon word, $u$ has no border by Proposition 2.2(iii),
and so $u$ cannot arise as a prefix power of a word $v\in n_{i}$.
Since $u$ is a prefix power of some word in some necklace of $M$,
it follows that $M$ contains a necklace of cardinal at least $n$.
$\qed$

The bound of $n$ in Lemma 3.3 is tight: see Theorem 3.8 below. It
follows from Lemma 3.3 that the length $l$ of the rows of the table
$T=T(M)$ is at least $n$. Consider the sub-table consisting of the
first $n$ columns of $T$. Since $M$ is an $n$ span de Bruijn set,
the rows of this sub-table form the dictionary of $A^{n}$. Each $u\in A^{n-1}$
is the prefix of $k$ successive rows of $T$ and if two of these
rows ended with the same letter $a\in A$, then the images of these
two rows under $\Pi$ would both begin with $au$, from which it would
follow that $au\in A^{n}$ would be a prefix of a power of two distinct
words of the necklaces of $M$, contrary to $M$ being a de Bruijn
set of span $n$. It follows that the final column of $T$ is a product
of $k^{n-1}$ members (possibly with repetitions) taken from the set
$G=\{i_{1}i_{2}\cdots i_{k}:\,\{i_{1},i_{2},\cdots,i_{k}\}=[k]\}$
(that is, $G$ consists of all $k!$ products of distinct members
of $A$). These observations establish the forward implication in
the following result.

\textbf{Theorem 3.4} The set of all BW transforms of de Bruijn sets
$M$ of span $n$ over a $k$-letter alphabet is $\Gamma_{k,n}=G^{k^{n-1}}$.

\textbf{Examples 3.5 }Let $k=2,\,n=4$. We may write $A=\{a<b\}$
so that $G=\{\alpha,\beta\}$ where $\alpha=ab,$ $\beta=ba$. Take
$v=\beta^{4}\alpha\beta^{3}\in G^{k^{n-1}}=G^{8}$. The standard permutation
$\pi(v)$ is the transitive cycle 
\[
\pi(v)=(0\,1\,3\,7\,15\,14\,12\,9\,2\,5\,11\,6\,13\,10\,4\,8),
\]
yielding the span $4$ Lyndon de Bruijn word $w=aaaa\,bbbb\,aaba\,bbab$.
As a second example take $v=\beta\alpha^{2}\beta^{2}\alpha^{2}\beta$
so that 
\[
\pi(v)=(0\,1\,2\,4\,9\,3\,7\,15\,14\,13\,11\,6\,12\,8)(5\,10);
\]
the corresponding set of Lyndon words is $\{aaaabaabbbbabb,\,ab\}$,
the cyclic $4$-factors of which are all the $2^{4}=16$ words of
$A^{4}$ with $\{abab,baba\}$ arising from the necklace defined by
the Lyndon word $ab$. 

We prove the reverse implication in Theorem 3.4 via two lemmas. 

\textbf{Lemma 3.6 }Let $v\in\Gamma_{k,n}$. Then $\pi(v)=\pi=\cup_{i=0}^{i=k-1}\pi_{i}$,
a union of $k$ order preserving partial mappings with dom $\pi_{i}=\{x=\varepsilon_{1}\varepsilon_{2}\cdots\varepsilon_{n}\in[k^{n}]:\varepsilon_{1}=i\}$
for $0\leq i\leq k-1.$ The sets ran $\pi_{i}$ also partition $[k^{n}]$
and each range set is itself a transversal of the partition of $[k^{n}]$
into the successive intervals of length $k$ which are: 
\begin{equation}
[jk,(j+1)k-1],\,0\leq j\leq k^{n-1}-1
\end{equation}

\emph{Proof} The description of the sets dom $\pi_{i}$ follows from
the fact that $|v|_{i}$ is the same value, $k^{n-1}$, for each $i$
$(0\leq i\leq k-1)$ and the sets ran $\pi_{i}$ always partition
the base set as $\pi$ is a permutation. The claim as regard transversals
follows as each $v\in\Gamma_{k,n}$ is a product of words from $G$.
$\qed$

For any $x\in[k^{n}]$ and integer $m\geq0$ there is a unique product
$p=p_{x,m}=\pi_{\varepsilon_{1}}\pi_{\varepsilon_{2}}\cdots\pi_{\varepsilon_{m}}$
with each $\varepsilon_{i}\in[k]$, such that $x\cdot p$ is defined.
The product $p_{x,m}$ can therefore be identified with $\varepsilon_{1}\varepsilon_{2}\cdots\varepsilon_{m}$,
which we shall call the\emph{ $m$-string }of $x$. 

\textbf{Lemma 3.7 }Let $e=\varepsilon_{1}\varepsilon_{2}\cdots\varepsilon_{m}$
be an $m$-digit $k$-ary expression $(1\leq m\leq n)$. Then for
any $x\in[k^{n}]$ whose $n$-digit $k$-ary representation has $e$
as a prefix, the $k$-ary $m$-string of $x$ is $e$ in the standard
permutation $\pi(v),$ for every $v\in\Gamma_{k,n}$. Moreover, the
domain of the partial mapping $p_{e}=\pi_{\varepsilon_{1}}\pi_{\varepsilon_{2}}\cdots\pi_{\varepsilon_{m}}$
is the interval of all $x$, the $k$-ary representation of which
begins with $e$. In particular, dom $p_{x}=\{x\}$. 

\emph{Proof} By Lemma 3.6, $x\cdot\pi_{i}$ is defined if and only
if $F(x)=i$ and so the claim holds if $m=1$. We shall now verify
that $x\cdot\pi_{\varepsilon_{1}}$ has the $k$-ary form $\varepsilon_{2}\varepsilon_{3}\cdots\varepsilon_{n}\varepsilon_{1}'$,
$(\varepsilon_{1}'\in[k])$, from which the result follows by repeated
application of this fact. Now since $\pi_{\varepsilon_{1}}$ is order-preserving,
it follows from Lemma 3.6 that we may identify the interval of (5)
in which $x\cdot\pi_{\varepsilon_{1}}$ lies by putting $j=\varepsilon_{2}\varepsilon_{3}\cdots\varepsilon_{n}$,
giving: 
\[
[(\varepsilon_{2}\varepsilon_{3}\cdots\varepsilon_{n})k,\,(\varepsilon_{2}\varepsilon_{3}\cdots\varepsilon_{n}+1)k-1]=[\varepsilon_{2}\varepsilon_{3}\cdots\varepsilon_{n}0,\,\varepsilon_{2}\varepsilon_{3}\cdots\varepsilon_{n}0+(k-1)]
\]
and so $x\cdot\pi_{\varepsilon_{1}}=\varepsilon_{2}\varepsilon_{3}\cdots\varepsilon_{n}\varepsilon_{1}'$
, as required. By what we have just proved and the uniqueness of the
products $p_{x,m}$, the integer $x\in$ dom $p_{e}$ if and only
if $e$ is a prefix of $x$, whence the final claim follows. $\qed$ 

\emph{Proof of Theorem 3.4}\textbf{. }The forward implication was
proved in the preamble to the theorem so consider the converse. For
$v\in\Gamma_{k,n}$ consider $M=BW^{-1}(v)$. By Lemma 3.7, for any
$x\in[k^{n}]$, $u=x$ is the unique word $u\in A^{n}$ such that
$x\cdot\pi_{u}$ is defined. Since some members $x\in[k^{n}]$ such
as $x=1$, are primitive, the table $T(M)=T(v)$ has at least $n$
columns. It follows that the prefix of length $n$ of the row $x$
of $T(v)$ is $x$ and so the sub-table of the first $n$ columns
of $T(v)$ has as its rows the members of $[k^{n}]$ written in numerical
order. In particular $x$ occurs among the $k^{n}$ factors of length
$n$ that can be read from the $k^{n}$ words of the necklaces of
$M$, and so each such $x$ must occur exactly once and therefore
$M$ is a de Bruijn set of span $n$. $\qed$

We next look at the special case where $v$ is a power of $\alpha=1\,2\cdots\,(k-1)$. 

\textbf{Theorem 3.8 }Let $v=\alpha^{k^{n-1}}$, let $M=BW^{-1}(v)$
and let $T=T(v)=T(M)$. Then the rows of $T$ are simply the list
of numbers $[k^{n}]$. Moreover $BW^{-1}(v)$ is the set of necklaces
of Lyndon words of length dividing $n$. The Lyndon words of the roots
of the necklaces of $M$ occur in the rows of $T$ in lexicographic
order. 

\emph{Proof} As in the proof of Theorem 3.4, we see that the sub-table
of the first $n$ columns of $T$ simply lists the numbers of $[k^{n}]$.
However since $v=\alpha^{k^{n-1}}$, for $x=\varepsilon_{1}\varepsilon_{2}\cdots\varepsilon_{n}$,
$x\cdot\pi_{\varepsilon_{1}}$ is the $\varepsilon_{1}$th member
$(0\leq\varepsilon_{1}\leq k-1)$ of the specified interval in (5),
that is to say, $x\pi_{\varepsilon_{1}}=\varepsilon_{2}\varepsilon_{3}\cdots\varepsilon_{n}\varepsilon_{1};$
by repetition of this observation we infer that for any $x=\varepsilon_{1}\varepsilon_{2}\cdots\varepsilon_{n}$,
the sequence $x,x\cdot\pi,x\cdot\pi^{2},\cdots x\cdot\pi^{n-1}$(where
$\pi=\pi(v)$) is the cyclic sequence under $\Pi$ of $x=\varepsilon_{1}\varepsilon_{2}\cdots\varepsilon_{n}.$
Since $x\cdot\pi^{n}=x$, it follows that the cardinal of the corresponding
necklace is a divisor of $n$; in particular $l$, the lcm of the
length of the roots of words of the rows is $n$, so that $T$ is
simply the table of $[k^{n}]$. The least member of each necklace
is by definition a Lyndon word. Every Lyndon word $w$ of length dividing
$n$ has a power which is some word $x\in[k^{n}]$ and so $w$ occurs
as a Lyndon word of some necklace in $BW^{-1}(v)$. The Lyndon roots
of the words of $T$ occur in lexicographic order by Theorem 1.2.13.
$\qed$

\textbf{Example 3.9 }Let us take $k=2,\,n=5$, and again reverting
to the alphabet $A=\{a<b\}$, we have $\alpha=ab$ and $v=\alpha^{2^{4}}=(ab)^{16}$.
Then 
\[
\pi(v)=(0)(1\,2\,4\,8\,16)(3\,6\,12\,24\,17)(5\,10\,20\,9\,18)(7\,14\,28\,25\,19)
\]
\[
(11\,22\,13\,26\,21)(15\,30\,29\,27\,23)(31).
\]
Expressed as a concatenation of Lyndon words of the corresponding
necklaces we obtain: 
\[
BW^{-1}(v)=a\cdot aaaab\cdot aaabb\cdot aabab\cdot aabbb\cdot ababb\cdot abbbb\cdot b.
\]
This is indeed the first de Bruijn word of span $5$ in the lexicographic
order. That this is always the case is a well-known theorem of Frederickson
and Maiorana. (See also {[}10{]} for an alternative proof.)

\textbf{Theorem 3.10} {[}4{]} For a given $n$, the lexicographic
concatenation of all Lyndon words of length dividing $n$ is the de
Bruijn word of span $n$ that lies first in the lexicographic order.

\textbf{Corollary 3.11 }Taken in ascending order of their Lyndon words,
the concatenation of the Lyndon words of the necklaces of\textbf{
}$BW^{-1}(\alpha^{k^{n-1}})$ is the first de Bruijn word of span
$n$ in the lexicographic order.

\textbf{Example 3.12 }Let us take $k=n=3$ so that $\alpha=abc$ say
and calculate $BW^{-1}(v)$ where $v=\alpha^{k^{n-1}}=(abc)^{9}.$
We find that 
\[
\pi_{v}=(0)(1\,3\,9)(2\,6\,18)(4\,12\,10)(5\,15\,19)(7\,21\,11)
\]
\[
(8\,24\,20)(13)(14\,16\,22)(17\,25\,23)(26);
\]
and so the least de Bruijn word that contains all words of length
$3$ over the alphabet $A=\{a<b<c\}$ as its set of cyclic factors
is the following concatentation of Lyndon words of lengths 1 or 3
over $A$:
\[
BW^{-1}(\alpha^{9})=a\cdot aab\cdot aac\cdot abb\cdot abc\cdot acb\cdot acc\cdot b\cdot bbc\cdot bcc\cdot c.
\]

\section{Maximum number of distinct factors of a word }

As an application of de Bruijn words we derive the functional form
for the maximum number of distinct factors in $A^{+}$of a word of
length $n$ over a fixed finite alphabet $A$. The upper bound in
our result comes from observing that long words must have repeated
short factors while the proof for the lower bound relies on the fact
that factors of de Bruijn words have no repeats of their long factors.
The topic of the number of subwords\emph{ }of a word has been extensively
investigated: for example see Section 6.3 of {[}7{]}. 

Consider the finite alphabet $A=A_{k}=\{a_{1},a_{2},\cdots,a_{k}\}$.
The set $A^{\leq m}=\{w:w\in A^{+}\,\mbox{ and \ensuremath{|w|\leq m\}}}$.
The number of distinct factors of $w$ will be denoted by $f_{w}$. 

\textbf{Lemma 4.1} With repeats, the number of factors in $A^{+}$of
$w\in A^{n}$ $(n\geq1)$is $\frac{1}{2}n(n+1)$.

\emph{Proof} A factor of $w$ is determined by the choice of two distinct
positions with each position occurring either between letters or at
either end of $w$. There are $\binom{n+1}{2}=\frac{1}{2}n(n+1)$
such pairs. $\qed$

\textbf{Corollary 4.2} For $w\in A^{n}\,(n\geq1)$ we have $n\leq f_{w}\leq\frac{1}{2}n(n+1)$.
Moreover, the lower bound is obtained if and only if $|c(w)|=1$ and
the upper bound is attained if and only if $n\leq k$.

\emph{Proof} The upper bound for $f_{w}$ comes from Lemma 4.1. Since
any word $w\in A^{n}$ has $n$ distinct prefixes it follows that
$n\leq f_{w}$ always holds. If $|c(w)|=1$, then $w=a^{n}$ for some
$a\in A$ and the set of factors of $w$ is $\{a^{t}:1\leq t\leq n\}$
and is of cardinal $n.$ On the other hand if $|c(w)|\geq2$ then,
in addition to its $n$ prefixes, $w$ also has the factor $b\in A$
where $b\neq F(w)$ so that $n<f_{w}$. Next suppose that $n\leq k$.
Put $w=a_{1}a_{2}\cdots a_{n}$; no two factors of $w$ have the same
content so the factors of $w$ are pairwise distinct, showing that
the upper bound in the statement is attained in this case. For all
remaining cases we have $2\leq k<n$ in which instance $w$ has two
identical $1$-factors and so $f_{w}<\frac{1}{2}n(n+1)$. $\qed$

In light of Corollary 4.2 we shall henceforth assume that $2\leq k<n$. 

\textbf{Definition 4.3} Let $f(n)=$ max$\{f_{w}:w\in A^{n}\}$.

\textbf{Theorem 4.4} $\frac{1}{2}n^{2}-f(n)=O(n\log n)$.

\emph{Proof} For $1\leq r\leq n$, a word $w\in A^{n}$ has $n-r+1$
(not necessarily distinct) $r$-factors and $|A^{r}|=k^{r}$. Hence
there are at least $n-r+1-k^{r}$ repeated $r$-factors in $w$. Let
$t$ be the greatest value of $r$ such that $r+k^{r}\leq n$, noting
that $1\leq t$. The total number of repeated factors in $w$ is then
at least:
\begin{equation}
\sum_{r=1}^{t}(n-r+1-k^{r})=(n+1)t-\frac{1}{2}t(t+1)-k\frac{k^{t}-1}{k-1}
\end{equation}
Now since $t+k^{t}\leq n<t+1+k^{t+1}$ we have $n<2k^{t+1}$; by taking
logarithms to the base $k$ we obtain $t<\log_{k}n<(1+\log_{k}2)+t$
so that $t=O(\log n).$ Moreover, $k^{t}=O(n)$ whence it follows
that 
\begin{equation}
f(n)\leq\frac{1}{2}n(n+1)-(n+1)O(\log n)+\frac{1}{2}(O(\log n))^{2}+O(n)=\frac{1}{2}n^{2}-O(n\log n)
\end{equation}

Conversely, given $n$, let $m\geq1$ be determined by the inequalities
$k^{m-1}<n\leq k^{m}$. Take $w\in A^{n}$ to be a factor of a de
Bruijn word $d=d_{m}$ of span $m$ over $A$. For any positive integer
$p\leq n$ there are $n-p+1$ factors of length $p$ in $w$. Moreover
if $m\leq p$, these factors are pairwise distinct as the members
of the set of prefixes of length $m$ of these factors are pairwise
distinct since $d$ is a de Bruijn word of index $m$. Hence 
\[
f_{w}\geq1+2+\cdots+(n-m+1)
\]
\begin{equation}
\Rightarrow f_{w}\geq\frac{1}{2}(n-m+1)(n-m+2)=\frac{1}{2}n^{2}-nm+\frac{1}{2}(3(n-m)+m^{2})+1>\frac{1}{2}n^{2}-nm
\end{equation}
Now we have $k^{m-1}<n\leq k^{m}$, whence $m=O(\log n)$ and so (8)
yields: 
\begin{equation}
f(n)\geq f_{w}\geq\frac{1}{2}n^{2}-O(n\log n)
\end{equation}
Combining (7) and (9) we conclude that $\frac{1}{2}n^{2}-f(n)=O(n\log n).$
$\qed$

\textbf{Acknowledgement} I would like to thank the referees for their
constructive suggestions and Alexei Vernitski for pointing out the
connection between the standard permutation and syntactic semigroups.


\begin{thebibliography}{References}
\bibitem{key-1} van Aardenne-Ehrenfest T. and de Bruijn N.G., \emph{Circuits
and trees in oriented linear graphs}, Simon Stevin \textbf{28}: 203-217
(1951).

\bibitem{key-5}Burrows M. and Wheeler D.J., \emph{A block sorting
data compression algorithm,} Technical Report, DIGITAL System Center,
(1994).

\bibitem{key-1}Crochemore M.J., Désarménien J. and Perrin D., \emph{A
note on the Burrows-Wheeler transformation}, Theoretical Computer
Science, \textbf{Vol. 332} Issue 1-3, 2005, 567-572.

\bibitem{key-7} Fredricksen H. and Maiorana J., \emph{Necklaces of
beads in k colors and k-ary de bruijn sequences}, Discrete Math. \textbf{23}
(1978), 207-210.

\bibitem{key-4}Gessel I.M. and Reutenauer C., \emph{Counting permutations
with given cycle structure and descent set, }J. Combin. Thry. Series
A, \textbf{64}, 189-215.

\bibitem{key-9}Higgins P.M., \emph{The semigroup of conjugates of
a word}, International Journal of Algebra and Computation, \textbf{Vol.
16}, No. 6 (2006), 1015-1029.

\bibitem{key-10}Lothaire M., `Combinatorics on Words', Cambridge
University Press, (2002). 

\bibitem{key-2}Mantaci S., Restivo A., Rosone G. and Sciortino M.,
\emph{An extension of the Burrows-Wheeler Transform and applications
to sequence comparison and data compression, }in `Combinatorial Pattern
Matching', Lecture Notes in Computer Science, (2005), Vol. 3537/2005,
178-189.

\bibitem{key-1}Mantaci S., Restivo A., Rosone G. and Sciortino M.,
\emph{An extension of the Burrows-Wheeler Transform}, Theor. Comput.
Sci. 387(3): 298-312 (2007). 

\bibitem{key-8}Moreno E., \emph{On the theorem of Fredricksen and
Maiorana about de Bruijn sequences} Adv. in Appl. Math. 33 (2), pp.413-415.

\bibitem{key-1}\emph{ }Pin J.E.,\emph{ Varieties of Formal Languages},
Plenum Publishing Corporation, New York, 1986. x + 138 pp. 
\end{thebibliography}
\end{document}